\documentclass[reqno]{amsart} 
\usepackage{amsfonts}
\usepackage{amssymb}
\usepackage{amsthm}

 \usepackage[dvips]{epsfig}  
 \usepackage{mathrsfs}
\usepackage{amssymb, amsmath}  
\usepackage[mathcal]{euscript}  
\usepackage{natbib}   
\numberwithin{equation}{section}   
   
\theoremstyle{plain}   
\newtheorem{thm}{Theorem}[section]   
\newtheorem{lem}[thm]{Lemma}   

\newtheorem{prop}[thm]{Proposition}
   
\theoremstyle{definition}

\theoremstyle{remark}   
\newtheorem{rem}[thm]{Remark}
\newcommand{\eqdis}{\stackrel{\lower0.2ex\hbox{$\scriptscriptstyle   
                    \mathrm{d}$}}{=}}   
\newcommand{\vague}{\stackrel{\lower0.2ex\hbox{$\scriptscriptstyle   
                    \it{v} $}}{\rightarrow}}   
\newcommand{\weak}{\stackrel{\lower0.2ex\hbox{$\scriptscriptstyle   
                    \it{w} $}}{\rightarrow}}   
\newcommand{\what}{\stackrel{\lower0.2ex\hbox{$\scriptscriptstyle   
                    \it{\hat{w}} $}}{\rightarrow}}   
\newcommand{\distr}{\stackrel{\lower0.2ex\hbox{$\scriptscriptstyle   
                    \it{d} $}}{\rightarrow}}   
\newcommand{\ntoinf}{\stackrel{\lower0.2ex\hbox{$\scriptscriptstyle   
                    \it{n\to\infty} $}}{\rightarrow}}   
\newcommand{\as}{\stackrel{\lower0.2ex\hbox{$\scriptscriptstyle   
                    \it{a.s.} $}}{\rightarrow}}   
\newcommand{\inprob}{\stackrel{\lower0.2ex\hbox{$\scriptscriptstyle   
      \Prob$}}{\rightarrow}}   
\newcommand{\Prob}{{P}}
\newcommand{\Var}{\operatorname{Var}}   
\newcommand{\Cov}{\operatorname{Cov}}   
\newcommand{\VaR}{\operatorname{VaR}}
\newcommand{\R}{\mathbf{R}}   
\newcommand{\vep}{\varepsilon}

 \makeatletter
 \AtBeginDocument{%
  \def\@serieslogo{%
  \vbox to\headheight{%
  \parindent\z@ \fontsize{6}{7\p@}\selectfont
  September 11, 2009\endgraf
  \vss}}}
 \makeatother

\begin{document}

\title[Efficient calculation of risk measures]{Efficient calculation
  of risk measures by importance sampling -- the heavy tailed case}


\author[H.~Hult]{Henrik Hult}  
\address[H.~Hult]{Department of Mathematics, KTH, 100 44 Stockholm, Sweden} 
\email{{hult@kth.se}}  
\author[J.~Svensson]{Jens Svensson}  
\address[J.~Svensson]{Department of Mathematics, KTH, 100 44
  Stockholm, Sweden}    
\email{{jenssve@kth.se}}


\begin{abstract}   
Computation of extreme quantiles and tail-based risk measures using
standard Monte Carlo simulation can be inefficient. A method to
speed up computations is provided by importance sampling. We show that
importance sampling algorithms, 
designed for efficient tail probability estimation, can significantly
improve Monte Carlo estimators of tail-based risk measures. In the
heavy-tailed setting, when the 
random variable of interest has 
a regularly varying distribution, we provide sufficient conditions for
the asymptotic relative error of importance sampling estimators of
risk measures, such as Value-at-Risk and expected shortfall, to be
small. The results are illustrated by some numerical examples. 
\end{abstract}   
   
   


\copyrightinfo{}{The authors}
 \maketitle   

\section{Introduction}

Risk measures are frequently used to quantify uncertainty in a
financial or actuarial context. Many risk measures, such as
Value-at-Risk and expected shortfall, depend on the tail
of the loss distribution. Exact formulas for computing risk
measures are only available for simple models and an alternative is to use
Monte Carlo simulation. However, standard Monte Carlo can be
inefficient when the function of interest depends on the
occurrence of rare events. A large number of samples may be
needed for accurate computation of extreme risk measures with standard
Monte Carlo, resulting in high computational cost. An alternative to
reduce the computational cost without 
loss of accuracy is provided by importance sampling. There is a vast
literature on the design of importance sampling algorithms for
computing rare event probabilities. However, for computing quantiles
and other risk measures the literature is not as
developed. \cite{GHS02} propose a method for  
efficient computation of quantiles of a heavy-tailed portfolio. Their
method is based on efficient algorithms designed for
rare event probabilities. The exceedance probability is computed for
a suggested initial value of the quantile. Then the quantile estimate
is updated depending on the computed probability and a search
algorithm is constructed to generate subsequent and more accurate
quantile estimates. We follow a more direct approach suggested
by \cite{G96} for computing quantiles and  extend it to handle
expected shortfall.

Let us first give a brief description of the problem. 
Let $X$ be a random variable with distribution $\mu$, distribution 
function (d.f.)~$F$, and continuous density $f$. Consider the problem
of computing its $p$th 
quantile, i.e. a number $\lambda_p$ such that 
$P(X\geq \lambda_p)\geq 1-p$, for some $p\in (0,1)$. If possible the
quantile is calculated by inverting  the d.f.
\begin{align*}
  \lambda_p= F^\leftarrow(p) := \inf \{x:1-F(x)\geq 1-p \}.
\end{align*}
When this is impossible an alternative is to use simulation. 
Computation of $F^\leftarrow(p)$ by standard Monte Carlo can be
implemented as follows. Generate $N$ independent copies
of $X$, denoted $X_1, \dots, X_N$. The empirical distribution function
(e.d.f.) of the sample is 
given by
\begin{align*}
  F_N(x)=\frac{1}{N}\sum_{i=1}^N I\{X_i\leq x \}
\end{align*}
and the quantile estimate is given by $F_N^\leftarrow(p)$. For
extreme quantiles, when $1-p$ is very small, standard Monte Carlo 
can be inefficient. Since only a small fraction of the sample will be
located in the tail, large samples are needed to obtain reliable estimates. 
A rough approach for quantifying the efficiency of Monte Carlo
estimates is to first consider a central limit theorem for
$F_N^\leftarrow(p)$. Suppose (this is true under suitable conditions)
\begin{align*}
  \sqrt{N}(F_N^\leftarrow(p)-F^\leftarrow(p)) \weak
  N\Big(0,\frac{p(1-p)}{f(F^\leftarrow(p))^2}\Big), 
\end{align*}
where $\weak$ denotes weak convergence. 
It is desirable to have the asymptotic standard deviation of roughly
the same size as the quantity $F^\leftarrow(p)$ we are estimating. For
standard Monte Carlo the asymptotic standard deviation is
\begin{align*}
  \frac{\sqrt{p(1-p)}}{f(F\leftarrow(p))} \approx 
  \frac{\sqrt{\overline{F}(F^\leftarrow(p))}}{f(F^\leftarrow(p))}
\end{align*}
which is typically much larger than $F^\leftarrow(p)$ for $p$ close to $1$. 
 
Consider, as an alternative to standard Monte Carlo, the method of
importance sampling. Then the sample $X_1,\ldots, 
X_N$ is generated from the sampling distribution $\nu$ and the
importance sampling tail e.d.f.~is given by 
\begin{align*}
  \overline{F}_{\nu,N}(x)=\frac{1}{N} \sum_{i=1}^N
  \frac{d\mu}{d\nu}(X_i) I\{X_i>x\}. 
\end{align*}
The quantile estimate is then given by
$(1-\overline{F}_{\nu,N})^\leftarrow(p) = \inf\{x:
\overline{F}_{\nu,N}(x) \leq 1-p\}$. The goal is to choose $\nu$ 
to get many samples in the tail of the original distribution and with
small Radon-Nikodym weights. Again, a rough evaluation of the
performance may be done by studying the limiting variance $\sigma_p^2$
in a central limit theorem of the form 
\begin{align*}
  \sqrt{N}((1-\overline{F}_{\nu,N})^\leftarrow(p)-F^\leftarrow(p)) \weak
  N\Big(0,\sigma_p^2\Big). 
\end{align*}
It turns out that the asymptotic properties (as $p \to 1$) of
$\sigma_p^2$ is closely related to the asymptotics of the second
moment of importance sampling algorithms designed for computing rare
event probabilities. This indicates that efficient algorithms for
computing rare event probabilities are indeed useful  for computing
quantiles. 

We use a similar approach to evaluate the performance of importance
sampling algorithms for computing expected shortfall. For  a random
variable $X$ with d.f.\ $F$ expected shortfall at level $p \in
(0,1)$ can be expressed as  
\begin{align*}
  \text{Expected shortfall}_p(X) = \frac{1}{1-p}\int_p^1
  F^\leftarrow(u)du =: \gamma_p(F^\leftarrow). 
\end{align*}
The standard Monte Carlo estimate based on a sample $X_1, \dots, X_N$
is given by $\gamma_p(F_N^\leftarrow)$ whereas the importance
sampling estimate is given by $\gamma_p((1-\overline{F}_{\nu,N})^\leftarrow)$.
A central limit theorem is derived for the expected shortfall estimate
based on importance sampling and the properties of the limiting
variance are studied as $p \to 1$. 

When evaluating the asymptotic variance for $p$ close to $1$ we
restrict attention to the heavy-tailed case. More precisely, it is
assumed that the original distribution has a regularly varying
tail. This is motivated by applications, e.g.\ finance and
insurance, where heavy-tailed data are frequently observed and
evaluation of risk measures is important for risk control.

For computation of rare event tail probabilities, the importance
sampling estimate of 
$P(X>\lambda)$ is given by $\hat
p_\lambda=\overline{F}_{\nu,N}(\lambda)$. Typically, the performance
of a rare event simulation algorithm 
is evaluated in terms of the relative error,  
\begin{align*}
  \textrm{Relative Error}=\frac{\sqrt{\Var(\hat p_\lambda)}}{p_\lambda}
\end{align*}
An algorithm is said to be \emph{asymptotically optimal} if the
relative error tends to 0 
as $\lambda\to\infty$. If the relative error remains bounded as
$\lambda\to\infty$, the 
algorithm is said to have \emph{bounded relative error}.
In the heavy-tailed setting there exist several algorithms for the case
where $X $ is given by the value at time $n$ of a random
walk. \cite{BJZ07} show that for 
such algorithms 
a necessary condition for them to achieve asymptotic optimality is that they are state-dependent.
\cite{DLW07} develop the first such algorithm, which almost
achieves asymptotic optimality for 
regularly varying distributions. 
\cite{BG08} propose a state-dependent algorithm with bounded
relative error for a more general class of 
heavy-tailed distributions.
\cite{BL08} consider the case where the number of steps of the
random walk and $\lambda$ tends to infinity simultaneously and 
develop an algorithm with 
bounded relative error.
\cite{HS09a} consider algorithms of the same kind as \cite{DLW07}, and show that they can
be made asymptotically optimal. 

The paper is organized as follows. In Section 2 
we review some standard results from empirical process theory. In Section 3 
we derive central limit theorems for empirical
quantiles for the empirical measures obtained from an importance
sampling algorithm. In Section 4 
we consider computation of risk measures for heavy-tailed (regularly varying) random
variables. Sufficient conditions for importance sampling algorithms
designed for rare event 
probability estimation to have small asymptotic variance are provided. 
Finally, in Section 5 
the procedure is illustrated for
computation of risk measures when the variable of interest is the
position of a finite random walk with regularly varying steps.

\section{Empirical processes}\label{sec:emp}

In this section we review some basic results from the theory of
empirical processes. We refer to \cite{VW96} for a thorough
introduction (see also \cite{CCH86}). 

Let $\{X_i\}_{i=1}^\infty$ be independent identically distributed random
variables with distribution $\mu$. Denote by $\mu_N$ the empirical
measure based on the first $N$ observations;
\begin{align*}
  \mu_N = \frac{1}{N}\sum_{i=1}^N \delta_{X_i},
\end{align*}
where $\delta_{x}$ is the Dirac measure at $x$. For a collection
$\mathcal{F}$ of real valued measurable functions, the empirical
measure induces a map from $\mathcal{F} \to \R$ by $f \mapsto
\mu_N(f) = \int f d\mu_N$. Assuming $\sup_{f \in \mathcal{F}}|f(x) -
\mu(f)| < \infty$ for each $x$, the empirical process $\xi_N$, given by
\begin{align*}
  \xi_N(f) = \sqrt{N}(\mu_N(f) - \mu(f)),
\end{align*}
can be viewed as a map into $l^\infty(\mathcal{F})$; the space of
bounded functions $\mathcal{F} \to \R$ equipped with the uniform metric. 
The collection $\mathcal{F}$ is called $\mu$-Donsker if 
\begin{align*}
  \xi_N \weak \xi,\quad \text{ in } l^\infty(\mathcal{F})
\end{align*}
and the limit is a tight Borel 
measurable element in $l^\infty(\mathcal{F})$.  

The classical result by Donsker states that the collection of indicator
functions $x \mapsto I\{x \leq t\}$ is $\mu$-Donsker with the limiting
process $B \circ \mu$, where $B$ is a Brownian bridge on $[0,1]$ (see
\cite{VW96}, pp. 81-82). In this paper we will be particularly
concerned with the collection $\mathcal{F}_{a,b}$ of indicator
functions $x \mapsto 
I\{a < x \leq t\}$ for $-\infty \leq a < t \leq b \leq \infty$, which
also is $\mu$-Donsker for any probability distribution $\mu$. To
simplify notation we will often write $\xi_N \weak \xi$ in
$l^\infty[a,b]$ for $\xi_N \weak \xi$ in $l^\infty(\mathcal{F}_{a,b})$.

To obtain convergence results for mappings of the empirical process
it is useful to apply the functional delta-method. 
Let $\mathbf{E}_1$ and $\mathbf{E}_2$ be two metric spaces. A mapping 
$\phi:\mathbf{E}_1 \to \mathbf{E}_2 $ is said to be {\it Hadamard differentiable
  at $\theta$ tangentially to $E_0\subset \mathbf{E}_1$} if there is a 
continuous mapping  $\phi_{\theta}'(h):\mathbf{E}_1 
\to \mathbf{E}_2 $ such that
\begin{align*}
  \frac{\phi(\theta + t_n h_n)-\phi(\theta)}{t_n}\to \phi_{\theta}'(h)
\end{align*}
for all sequences $t_n\to 0$ and $h_n\to h$, where $h\in E_0$.

\begin{thm}[Functional delta-method, c.f.~\cite{VW96}, Theorem 3.9.4]
  Let  $\phi:\mathbf{E}_1 \to \mathbf{E}_2 $ be Hadamard differentiable
  at $\theta$ tangentially to $E_0\subset \mathbf{E}_1$.
  Let $\{X_n\}_{n=1}^\infty$ be a sequence of random variables taking values in
  $\mathbf{E}_1$. Suppose $r_n(X_n-\theta)\weak X \in E_0$ for some
  sequence of constants $r_n\to\infty$.   
  Then $r_n(\phi(X_n)-\phi(\theta))\weak \phi_{\theta}'(X)$.
\end{thm}

For any c\`adl\`ag function $F: \R \to \R$, define the inverse map
$\phi_p$ by   
\begin{align*}
  \phi_p(F)&= F^\leftarrow(p) = \inf\{u:F(u)\geq p\}, \quad p \in (0,1). 
\end{align*}
The following result shows that the functional
delta-method implies the convergence of quantiles. 
\begin{prop}[c.f.~Lemma 3.9.23 and Example 3.9.24 in
  \cite{VW96}]\label{prop:quantile} 
  Let $\{X_i\}_{i=1}^\infty$ be a sequence of independent and identically
  distributed random variables with d.f.\ $F$. Suppose $F$ has a
  continuous density 
  $f>0$ with respect to the Lebesgue measure on the interval
  $[{F}^{\leftarrow}(p)-\epsilon,{F}^{\leftarrow}(q)+\epsilon]$, for
  $0<p<q<1$ and $\epsilon>0$.  
  Then   
  \begin{align*}
    \sqrt N ({F}_N^{\leftarrow}-{F}^{\leftarrow}) \weak \frac{B}{f \circ
      F^{\leftarrow}}, \quad \text{ in } l^\infty[p,q],
  \end{align*}
where the right-hand side refers to the random function $u \mapsto
\frac{B(u)}{f(F^\leftarrow(u))}$. 
\end{prop}

\section{Empirical processes and importance sampling} 
\label{sec:isemp}
The empirical measure resulting from a random sample of an
importance sampling algorithm with 
sampling distribution $\nu$ can be used to approximate important parts of the
original distribution. In our context $\nu$ is chosen to give a good
approximation of the extreme tail of the original distribution.  

Let $\{X_i\}_{i=1}^\infty$ be independent
identically distributed with distribution $\nu$. The empirical measure
with likelihood ratio weights and the corresponding tail empirical
distribution are written
\begin{align}
  \mu_{N,\nu} &= \frac{1}{N}\sum_{i=1}^N w(X_i) \delta_{X_i} =
  \frac{1}{N}\sum_{i=1}^N \frac{d\mu}{d\nu}(X_i) \delta_{X_i},
  \nonumber \\ 
   \overline{F}_{\nu,N}(t) &= \mu_{N,\nu}(I\{\cdot > t\}). \label{eq:tailedf}
\end{align}
Let $\mathcal{F}_a$ be the collection of indicator functions
$I\{\cdot > t\}$ with $t \geq a$. 
For the importance sampling estimators we are concerned with a central
limit theorem of the form 
\begin{align*}
  \sqrt{N} (\mu_{N, \nu} - \mu) \weak Z \quad \text{in } l^\infty(\mathcal{F}_a),
\end{align*}
which we write, with slight abuse of notation, as
\begin{align*}
  \sqrt N(\overline{F}_{\nu,N}-\overline{F})\weak Z, \quad \text{in }
  l^\infty[a,\infty]. 
\end{align*}
Note that  $\mu_{N,\nu}(f) = \nu_N(wf)$ and $\mu(f) = \nu(wf)$, where
$w = d\mu/d\nu$. Therefore
the central limit theorem can be stated by saying that the collection
$w\mathcal{F}_a = \{wf : f \in \mathcal{F}_a\}$ is $\nu$-Donsker.  By
the permanence properties of Donsker classes (see \cite{VW96}
Section 2.10) this follows when
$\mathcal{F}_a$ is $\nu$-Donsker and $E_\nu w(X)^2I\{X > a\} < 
\infty$. 

To identify the limiting process $Z$ we first need to calculate the
covariance function of the process $\overline{F}_{\nu,N}$.  
\begin{lem}\label{lem:covfcn}
  Let $\{X_i\}_{i=1}^\infty$ be independent and identically
  distributed with distribution $\nu$ with $\mu \ll \nu$ and $w =
  d\mu/d\nu$.  If $E 
  w(X_1)^2I\{X_1 > a\} < \infty$, for some 
  $a \geq - \infty$, then, for $y \geq x > a$,     
  \begin{align}\label{eq:rho}
    \varrho(x,y):=N \Cov(\overline{F}_{\nu,N}(x),\overline{F}_{\nu,N}(y)) 
    = E_{\nu}w(X_1)^2 I\{X_1>y\}
    - \overline{F}(x)\overline{F}(y). 
  \end{align}
\end{lem}

\begin{proof}
This is a direct calculation. Indeed,
  \begin{align*}
    \Cov&(\overline{F}_{\nu,N}(x),\overline{F}_{\nu,N}(y))\\
    &=E_{\nu}(\overline{F}_{\nu,N}(x) \overline{F}_{\nu,N}(y))-
    E_{\nu}(\overline{F}_{\nu,N}(x))E_{\nu}(\overline{F}_{\nu,N}(y))\\
    &= \frac{1}{N^2}\sum_{i,j}E_{\nu}( w(X_i)  
    w(X_j)I\{ X_i>x \} I\{ X_j>y \}) - \overline{F}(x)\overline{F}(y)\\
    &=\frac{1}{N} E_{\nu}( w(X_1)^2 I\{ X_1>x \}I\{ X_1>y \})\\
    &\quad+\frac{1}{N^2}\sum_{i\neq j}E_{\nu}( w(X_i)  w(X_j)I\{
    X_i>x \} I\{ X_j>y \}) 
    - \overline{F}(x)\overline{F}(y)\\
    &=\frac{1}{N} E_{\nu}( w(X_1)^2 I\{ X_1>y \})\\
    &\quad+\frac{N^2-N}{N^2}E_{\nu}( w(X_i)  w(X_j)I\{ X_i>x \} I\{ X_j>y \})
    - \overline{F}(x)\overline{F}(y)\\
    &=\frac{1}{N} E_{\nu}( w(X_1)^2 I\{ X_1>y
    \})+\frac{N^2-N}{N^2}\overline{F}(x)\overline{F}(y) -
    \overline{F}(x)\overline{F}(y)\\ 
    &=\frac{1}{N} E_{\nu}( w(X_1)^2 I\{ X_1>y \})
    -\frac{1}{N}\overline{F}(x)\overline{F}(y). 
  \end{align*}
\end{proof}
Note that if $w=1$, i.e.\ the sampling measure is the original
measure, then the covariance function  becomes 
$\overline{F}(y)-\overline{F}(y)\overline{F}(x)=\overline{F}(y)(1-\overline{F}(x))=F(x)(1-F(y)),
\, y>x$, 
which corresponds to a Brownian bridge evaluated at $F$.

Now we are ready to state the central limit theorem for the tail
empirical distribution of an importance sampling algorithm. 

\begin{prop}\label{prop:tedfconv}
  Let $Z$ be a centered Gaussian process with covariance function
  $\varrho$ given by \eqref{eq:rho}. If $E_\nu w(X_1)^2I\{X_1>a\} < \infty$ for
  some $a \geq -\infty$, then
  \begin{align*}
    \sqrt N(\overline{F}_{\nu,N}-\overline{F}) \weak Z, \quad \text{in }
    l^\infty[a,\infty]. 
  \end{align*}
\end{prop}
\begin{proof}
  We have already seen that $E_\nu w(X_1)^2I\{X_1>a\} <
  \infty$ implies that $\mathcal{F}_a$ is $\nu$-Donsker. Hence, we
  need only to identify the limiting process $Z$. Denote $\xi_N(x) = \sqrt 
  N(\overline{F}_{\nu,N}(x)-F(x))$.  
  By the multivariate central limit theorem the finite dimensional
  distributions converge; for any $x_1, \dots, x_k$ with $x_i > a$, 
  \begin{align*}
    (\xi_N(x_1), \ldots, \xi_N(x_k)) \weak N(0, \Sigma),
  \end{align*}
  where the entries of $\Sigma_{ij} = \varrho(x_i,x_j)$. This
  determines that the limiting process must be $Z$. 
\end{proof}

We proceed with the asymptotic normality of the quantile transform. 
The proof is very similar to that of Proposition \ref{prop:quantile} and
therefore omitted. 

\begin{prop}\label{prop:tedfinvconv}
  Let $Z$ be a centered Gaussian process with covariance function
  $\varrho$ in \eqref{eq:rho}.  Suppose $F$ has a continuous density 
  $f>0$ with respect to the Lebesgue measure on the interval
  $[{F}^{\leftarrow}(p)-\epsilon,{F}^{\leftarrow}(q)+\epsilon]$, for
  $0<p<q<1$ and $\epsilon>0$. $E_\nu
  w(X_1)^2I\{X_1>F^\leftarrow(p)-\vep\} < \infty$, then  
  \begin{align*}
    \sqrt N((1-\overline{F}_{\nu,N})^{\leftarrow}-{F}^{\leftarrow})\weak
    \frac{Z(F^{\leftarrow})}{f({F}^{\leftarrow})}, \quad \text{ in } l^\infty[p,q]. 
\end{align*}
\end{prop}

Next consider a central limit theorem for an importance sampling
estimate of expected shortfall. For a non-decreasing c\`adl\`ag
function $F^{\leftarrow}$ and $0 < p < 
1$ we use the notation
\begin{align*}
  \gamma_p(F^\leftarrow) &=\frac{1}{1-p}\int_p^1 F^{\leftarrow}(u)du. 
\end{align*}
Recall that expected shortfall at level $p$ for a random variable $X$
with d.f.\ $F$ is given by $\gamma_p(F^\leftarrow)$ and the importance
sampling estimate based on a sample $X_1, \dots, X_N$ with sampling
distribution $\nu$ is given by
$\gamma_p((1-\overline{F}_{\nu,N})^\leftarrow)$. 
\begin{prop}\label{prop:es}
  Assume the hypotheses of Proposition \ref{prop:tedfinvconv}. If, in
  addition,\\ 
  $\int_{F^\leftarrow(p)}^\infty\int_{F^\leftarrow(p)}^\infty
  \varrho(x,y)dxdy < \infty$ and $\varrho(x,x) =
  o([f(x)/\overline{F}(x)]^2)$, as $x \to \infty$, then 
  \begin{align*}
    \sqrt
    N(\gamma_p((1-\overline{F}_{\nu,N})^\leftarrow)-\gamma_p({F}^{\leftarrow}))
    \weak 
    \frac{1}{1-p}\int_p^1
    \frac{Z(F^{\leftarrow}(u))}{f({F}^{\leftarrow}(u))}du,
\end{align*}
as $N \to \infty$. 
\end{prop}
\begin{proof}
Let $q$ and $\vep$ be arbitrary with $p < q < 1$ and $\vep > 0$. Since
\begin{align}
  &
  P\Big(\Big|\sqrt{N}[\gamma_p((1-\overline{F}_{\nu,N})^\leftarrow)-
  \gamma_p(F^\leftarrow)]-
  \frac{1}{1-p}\int_p^1\frac{Z(F^\leftarrow(u))}{f(F^\leftarrow(u))}du\Big|
  > \vep\Big) \nonumber \\
  &\quad \leq
  P\Big(\Big|\frac{\sqrt{N}}{1-p} 
  \Big[\int_p^q(1-\overline{F}_{\nu,N})^\leftarrow(u)- 
  \int_p^q F^\leftarrow(u)du\Big]\nonumber  \\ & \qquad \qquad \qquad
  \qquad \qquad \qquad 
  \qquad \qquad - 
  \frac{1}{1-p}\int_p^q\frac{Z(F^\leftarrow(u))}{f(F^\leftarrow(u))}du\Big|
  > \vep/3\Big) \label{eq:empesthree.1}\\
  & \qquad + 
  P\Big(\sqrt{N}\Big|\frac{1}{1-p}\int_q^1
  (1-\overline{F}_{\nu,N})^{\leftarrow}(u)-{F}^{\leftarrow}(u)du
  \Big|> \vep/3\Big)  \label{eq:empesthree.2} \\
  & \qquad + 
  P\Big(\Big|\frac{1}{1-p}\int_q^1 \frac{Z(F^\leftarrow(u))}{f(F^\leftarrow(u))}du
  \Big|> \vep/3\Big) \label{eq:empesthree.3}
\end{align}
it is sufficient to show that each of the three terms converges to
$0$, as first $N \to \infty$ and then $q \to 1$.  

Consider first  \eqref{eq:empesthree.1}. Let $\gamma_{p,q}$ be the map defined by
\begin{align*}
  \gamma_{p,q}(H) &= \frac{1}{1-p}\int_p^q H(u)du,
\end{align*}
on the set $D_\gamma$ of all non-decreasing c\`adl\`ag functions
$H$. Since $\gamma_{p,q}$ is linear it is Hadamard differentiable on
$D_\gamma$ with derivative $\gamma'_{p,q}(h) = \gamma_{p,q}(h)$. 
In particular, Proposition \ref{prop:tedfinvconv} and the delta-method
imply that
  \begin{align*}
    \sqrt{N}(\gamma_{p,q}((1-\overline{F}_{\nu,N})^\leftarrow)-
    \gamma_{p,q}({F}^\leftarrow)) \weak  
    \frac{1}{1-p}\int_p^q
    \frac{Z(F^{\leftarrow}(u))}{f({F}^{\leftarrow}(u))}du.  
\end{align*}
This takes care of \eqref{eq:empesthree.1}. Next consider
\eqref{eq:empesthree.3}. By Chebyshev's inequality
\begin{align*}
  & P\Big(\Big|\frac{1}{1-p}\int_q^1 \frac{Z(F^\leftarrow(u))}{f(F^\leftarrow(u))}du
  \Big|> \vep/3\Big) \\ & \quad \leq
  \Big(\frac{3(1-p)}{\vep}\Big)^2 \Var\Big(\int_q^1
  \frac{Z(F^\leftarrow(u))}{f(F^\leftarrow(u))}du \Big)  \\ & \quad \leq
  \Big(\frac{3(1-p)}{\vep}\Big)^2 \int_q^1\int_q^1
  \frac{\varrho(F^\leftarrow(u),
    F^\leftarrow(v))}{f(F^\leftarrow(u))f(F^\leftarrow(v))}dudv  \\ &
  \quad =  
  \Big(\frac{3(1-p)}{\vep}\Big)^2
  \int_{F^\leftarrow(q)}^\infty\int_{F^\leftarrow(q)}^\infty \varrho(x,y) dxdy.
\end{align*}
Since the integral is finite, this converges to $0$ as $q \to 1$.

It remains to consider \eqref{eq:empesthree.2}.
First, write
\begin{align}
  \sqrt{N}&\Big|\frac{1}{1-p}\int_q^1
  (1-\overline{F}_{\nu,N})^{\leftarrow}(u)-{F}^{\leftarrow}(u)du\Big|
  \nonumber \\ &\leq
  \frac{\sqrt{N}}{1-p}
  \Big|\int_{(1-\overline{F}_{\nu,N})^{\leftarrow}(q)}^\infty
  \overline{F}_{\nu,N}(x) dx - \int_{F^{\leftarrow}(q)}^\infty
  \overline{F}(x)dx\Big| \nonumber \\ & \quad +
  \frac{\sqrt{N}(1-q)}{1-p}\Big|(1-\overline{F}_{\nu,N})^\leftarrow(q)-F^\leftarrow(q)\Big| 
  \nonumber \\  
  &\leq 
  \frac{\sqrt{N}}{1-p}
  \Big|\int_{F^{\leftarrow}(q)}^\infty \overline{F}_{\nu,N}(x) -
  \overline{F}(x)dx\Big| \label{eq:three.1} \\ & \quad +
  \frac{\sqrt{N}}{1-p}
  \Big|\int_{(1-\overline{F}_{\nu,N})^{\leftarrow}(q)}^{F^{\leftarrow}(q)}
  \overline{F}_{\nu,N}(x)dx\Big|I\{(1-\overline{F}_{\nu,N})^{\leftarrow}(q)
  \leq F^{\leftarrow}(q)\}  \label{eq:three.2.1} \\ &
  \quad +  \frac{\sqrt{N}}{1-p}
  \Big|\int_{F^{\leftarrow}(q)}^{(1-\overline{F}_{\nu,N})^{\leftarrow}(q)}
  \overline{F}_{\nu,N}(x)dx\Big| I\{(1-\overline{F}_{\nu,N})^{\leftarrow}(q)
  \geq F^{\leftarrow}(q)\} \label{eq:three.2.2} \\ & 
  \quad +
  \frac{\sqrt{N}(1-q)}{1-p}\Big|(1-F_{\nu,N})^\leftarrow(q)-F^\leftarrow(q)\Big|. 
  \label{eq:three.3} 
\end{align}
First consider \eqref{eq:three.1}. By Proposition \ref{prop:tedfconv}
and the delta method 
\begin{align*}
  \lim_{N\to\infty} & P\Big(\frac{\sqrt{N}}{1-p}
  \Big|\int_{F^{\leftarrow}(q)}^\infty \overline{F}_{\nu,N}(x) -
  \overline{F}(x)dx\Big| > \vep/12\Big) \\
  & = P\Big(\frac{1}{1-p}\Big|\int_{F^{-1}(q)}^\infty Z(x) dx\Big|
  > \vep/12\Big)\\
  & \leq \Big(\frac{12}{\vep (1-p)}\Big)^2
  \Var\Big(\int_{F^{\leftarrow}(q)}^\infty Z(x) dx \Big)\\
  & = \Big(\frac{12}{\vep (1-p)}\Big)^2
  \int_{F^{\leftarrow}(q)}^\infty\int_{F^{\leftarrow}(q)}^\infty \varrho(x,y) dxdy. 
\end{align*}
Since $F^{\leftarrow}(q) \to \infty$ as $q\to 1$ and the integral is
finite the expression in the last display can be made
arbitrarily small.  Next, consider \eqref{eq:three.2.1}. This term is
bounded from above by
\begin{align*}
  &\frac{\sqrt{N}}{1-p}\overline{F}_{\nu,N}\Big((1-\overline{F}_{\nu,N})^\leftarrow(q)\Big)\Big(F^\leftarrow(q) 
  - (1-\overline{F}_{\nu,N})^\leftarrow(q)\Big) \\ & \quad = \frac{1-q}{1-p}\sqrt{N}(F^\leftarrow(q)
  - (1-\overline{F}_{\nu,N})^\leftarrow(q)),
\end{align*}
where we have used that
$\overline{F}_{\nu,N}((1-\overline{F}_{\nu,N})^\leftarrow(q))
\geq 1-q$.  
By Proposition \ref{prop:tedfinvconv}, 
\begin{align*}
  &\lim_{N\to \infty}P\Big(\frac{1-q}{1-p}\sqrt{N}|F^\leftarrow(q)
  - (1-\overline{F}_{\nu,N})^\leftarrow(q)| > \vep/12\Big) \\ & \quad =
  P\Big(\frac{1-q}{1-p}\Big|\frac{Z(F^\leftarrow(q))}{f(F^\leftarrow(q))}\Big|
  > \vep/12\Big) \\ & \quad \leq
  \Big(\frac{12}{(1-p)\vep}\Big)^2(1-q)^2
  \frac{\varrho(F^\leftarrow(q),F^\leftarrow(q))}{f(F^\leftarrow(q))^2}\\ 
  & \quad =
  \Big(\frac{12}{(1-p)\vep}\Big)^2\overline{F}(F^\leftarrow(q))^2\frac{\varrho(F^\leftarrow(q),F^\leftarrow(q))}{f(F^\leftarrow(q))^2}.   
\end{align*}
This converges to $0$ as $q \to 1$ since $\varrho(x,x) =
o([f(x)/\overline{F}(x)]^2)$. 

Similarly, \eqref{eq:three.2.2} is
bounded from above by
\begin{align*}
  &\frac{\sqrt{N}}{1-p}\overline{F}_{\nu,N}(F^\leftarrow(q))\Big((1-\overline{F}_{\nu,N})^\leftarrow(q)
  - F^\leftarrow(q) \Big).
\end{align*}
This can be treated just like the previous term since, by Proposition
\ref{prop:tedfinvconv},  
\begin{align*}
  &\lim_{N\to \infty}P\Big(\frac{\overline{F}_{\nu,N}(F^\leftarrow(q))}{1-p}\sqrt{N}|(1-\overline{F}_{\nu,N})^\leftarrow(q)-F^\leftarrow(q)|
  > \vep/12\Big) \\ & \quad =
  P\Big(\frac{1-q}{1-p}\Big|\frac{Z(F^\leftarrow(q))}{f(F^\leftarrow(q))}\Big|
  > \vep/12\Big).
\end{align*}
Finally, \eqref{eq:three.3} can be treated the same way since, by
Proposition \ref{prop:tedfinvconv}, 
\begin{align*}
  & \lim_{N\to
    \infty}\Prob\Big(\frac{\sqrt{N}(1-q)}{1-p}\Big|(1-F_{\nu,N})^\leftarrow(q)-F^\leftarrow(q)\Big|  
  > \vep/12\Big) \\ & \quad =
  \Prob\Big(\frac{1-q}{1-p}\Big|\frac{Z(F^{\leftarrow}(q))}{f({F}^{\leftarrow}(q))}\Big| 
  > \vep/12\Big).
\end{align*}
This
completes the proof.
\end{proof}

\section{Efficient calculation of risk measures in the heavy-tailed setting}
\label{sec:rm}
In the previous section we established central limit theorems for
importance sampling estimates of Value-at-Risk (i.e.\ quantiles) and
expected shortfall. In this section we study the limiting variance of
the central limit theorems as a function of $p$ when $p$ is close to
$1$. The main requirement is that the asymptotic standard deviation
coming from the central limit theorem is
roughly of the same size as the quantity we are trying to compute,
when $p$ is close to $1$. We only consider the case
when the original distribution is heavy-tailed, in the sense
that $\overline{F}$ is regularly varying. 

There are three main assumptions in this section.

\begin{itemize}
\item We assume that the original distribution of interest
  has a regularly varying tail. That is, there exists $\alpha >
  0$ such that 
  \begin{align}\label{eq:as.regvar}
    \lim_{t \to \infty}\frac{\overline{F}(tx)}{\overline{F}(t)} =
    x^{-\alpha},\quad x > 0. 
  \end{align}
\item We assume that there is an available explicit asymptotic
  approximation for $\overline{F}(x)$. More precisely, we know a
  non-increasing function $U$ such that $U \sim \overline{F}$, i.e. 
  \begin{align}\label{eq:as.approx}
    \lim_{t \to \infty}\frac{U(t)}{\overline{F}(t)} = 1. 
  \end{align}
\item We assume that we can construct sampling measures $\nu_\lambda$
  with bounded relative error for computing rare event
  probabilities of the type  $\overline{F}(\lambda) = P(X > \lambda)$; i.e.
  \begin{align}\label{eq:as.bre}
    \limsup_{\lambda \to \infty}\frac{E_{\nu_\lambda}w_\lambda(X_1)^2
      I\{X_1> \lambda\}}{\overline{F}(\lambda)^2} < \infty.
  \end{align}
\end{itemize}

\subsection{Computation of quantiles -- Value-at-Risk}

For a random variable $X$ with d.f.\ $F$ the Value-at-Risk at level
$p$ is defined as the $p$th quantile; $\VaR_p(X) = F^\leftarrow(p)$. 
Given $p \in (0,1)$ close to $1$, the importance sampling estimate
based on independent and identically distributed samples with sampling
distribution $\nu$ is given by $(1-\overline{F}_{\nu,N})^\leftarrow
(p)$. Then Proposition \ref{prop:tedfinvconv} and Lemma \ref{lem:covfcn} 
determines the asymptotic variance as
\begin{align*}
  \Var\Big(\frac{Z(F^{\leftarrow}(p))}{f(F^\leftarrow(p))}\Big) &=
  \frac{\varrho(F^{\leftarrow}(p),F^\leftarrow(p))}{f(F^\leftarrow(p))^2}
  \\ &= \frac{E_{\nu} w(X_1)^2 I\{X_1 > F^{\leftarrow}(p)\}
    -\overline{F}(F^\leftarrow(p))^2}{f(F^\leftarrow(p))^2}
  \\ &= \frac{E_{\nu}w(X_1)^2 I\{X_1 > F^{\leftarrow}(p)\}
    -(1-p)^2}{f(F^\leftarrow(p))^2} := \sigma^2_p.
\end{align*}

To control the asymptotic variance it seems like a good choice to use
an efficient rare event simulation algorithm designed for efficient
computation of $P(X > F^\leftarrow(p)$. That is, with sampling
distribution $\nu_{F^{\leftarrow}(p)}$. This is of course
impossible since $F^{\leftarrow}(p)$ is unknown. However, the
asymptotic approximation $U$ may be helpful. Note that since $U$ is
monotone it has an inverse $U^\leftarrow$ and by regular variation
$(1-U)^\leftarrow \sim F^\leftarrow$ as $p \to 1$. Thus, it seems
reasonable to use $\nu_{u_p}$ where $u_p = (1-U)^\leftarrow(p)$. This
is justified by the next result.

\begin{prop} \label{prop:varas}
  Suppose \eqref{eq:as.regvar}-\eqref{eq:as.bre} hold. If there
  exists $c_0 < 1$ such that 
  \begin{align}\label{eq:as.rnscale}
    \limsup_{\lambda \to \infty}\frac{E_{\nu_\lambda}w_\lambda(X_1)^2
      I\{X_1 > c_0\lambda\}}{\overline{F}(\lambda)^2} < \infty,
  \end{align}
  then the
  sampling measures $\nu_{u_p}$ satisfy
  \begin{align*}
    \limsup_{p \to 1}\frac{\sigma_p^2}{F^\leftarrow(p)^2} < \infty. 
  \end{align*}
\end{prop}

\begin{proof}
First note that \eqref{eq:as.rnscale} implies that 
\begin{align*}
   \limsup_{\lambda \to \infty}\frac{E_{\nu_\lambda}w_\lambda(X_1)^2
      I\{X_1 > c\lambda\}}{\overline{F}(\lambda)^2} < \infty,
\end{align*}
for any $c \geq c_0$. Take $\vep \in (0,1-c_0)$. Then there exists $p_0$
such that $\frac{F^\leftarrow(p)}{u_p} > 1-\vep$, for each
$p \geq p_0$. In particular, 
\begin{align}
  &\limsup_{p \to 1} \frac{E_{\nu_{u_p}}w_{u_p}(X_1)^2
    I\{X_1 > F^{\leftarrow}(p)\}}{\overline{F}(u_p)^2} \nonumber
  \\ & \quad \leq  \limsup_{p \to 1} \frac{E_{\nu_{u_p}}w_{u_p}(X_1)^2
    I\{X_1 >
    (1-\vep)u_p\}}{\overline{F}(u_p)^2}
  < \infty \label{eq:as.rnscale1}
\end{align}
Since $u_p \sim F^\leftarrow(p)$ and $x f(x) \sim \alpha
\overline{F}(x)$, by Karamata's theorem, it follows that
\begin{align}
  \lim_{p \to 1} \frac{\overline{F}(U^{\leftarrow}(p))^2}{F^\leftarrow(p)^2
    f(F^\leftarrow(p))^2} = \lim_{p \to 1}
  \frac{\overline{F}(F^{\leftarrow}(p))^2}{F^\leftarrow(p)^2 
    f(F^\leftarrow(p))^2} &=
  \frac{1}{\alpha^2},  \label{eq:as.rnscale2} \\
  \lim_{p \to 1} \frac{(1-p)^2}{F^\leftarrow(p)^2
    f(F^\leftarrow(p))^2} =  \lim_{p \to 1}
  \frac{\overline{F}(F^{\leftarrow}(p))^2}{F^\leftarrow(p)^2 
    f(F^\leftarrow(p))^2} &= \frac{1}{\alpha^2}. \label{eq:as.rnscale3}
\end{align}
By \eqref{eq:as.rnscale1}, \eqref{eq:as.rnscale2}, and
\eqref{eq:as.rnscale3} it follows that
\begin{align*}
  \limsup_{p \to 1}\frac{\sigma^2_p}{F^\leftarrow(p)^2} &=
   \limsup_{p \to 1} \frac{E_{\nu_{u_p}}w_{u_p}(X_1)^2
    I\{X_1 > F^{\leftarrow}(p)\} - (1-p)^2}{F^\leftarrow(p)^2
    f(F^\leftarrow(p))^2} \\ 
  &=  \limsup_{p \to 1} \frac{E_{\nu_{u_p}}w_{u_p}(X_1)^2
    I\{X_1 > F^{\leftarrow}(p)\}}{\overline{F}(U^{\leftarrow}(p))^2}
  \frac{\overline{F}(U^{\leftarrow}(p))^2}{F^\leftarrow(p)^2
    f(F^\leftarrow(p))^2} \\
  & \quad - \frac{(1-p)^2}{F^\leftarrow(p)^2
    f(F^\leftarrow(p))^2} < \infty. 
\end{align*}
\end{proof}

Under somewhat stronger assumptions it is possible to reach a more
explicit asymptotic bound for $\frac{\sigma_p^2}{F^\leftarrow(p)^2}$.  

\begin{prop}\label{prop2:rnscale}
  Suppose \eqref{eq:as.regvar} and \eqref{eq:as.approx} hold. Suppose
  additionally that there exist $c_0 < 1$ and a function $\varphi$,
  continuous at $1$, such that, for $c \geq c_0$, 
  \begin{align}\label{eq:prop2}
    \lim_{\lambda \to \infty}\frac{E_{\nu_\lambda}w_\lambda(X_1)^2
      I\{X_1 > c\lambda\}}{\overline{F}(\lambda)^2} \leq \varphi(c).
  \end{align} 
  Then the sampling measures $\nu_{u_p}$ satisfy
  \begin{align*}
    \lim_{p \to 1}\frac{\sigma_p^2}{F^\leftarrow(p)^2} \leq
    \frac{\varphi(1)-1}{\alpha^2}. 
  \end{align*}
\end{prop}
\begin{rem}
  If the asymptotic quantile approximation based on $U$ always
  underestimates the true quantile, i.e.\ $u_p \leq
  F^\leftarrow(p)$ for each $p$, then one can take $c_0 = 
  1$ in Proposition \ref{prop:varas} and Proposition \ref{prop2:rnscale}.  
\end{rem}

\begin{proof}
First note that $(1-U)^\leftarrow \sim F^\leftarrow$. Take $\vep \in
(0,1-c_0)$. Then there exists $p_0$ 
such that $\frac{F^\leftarrow(p)}{u_p} > 1-\vep$, for each
$p \geq p_0$. In particular,
\begin{align*}
   \limsup_{p \to 1}\frac{E_{\nu_{u_p}}w_{u_p}(X_1)^2
      I\{X_1 > F^\leftarrow(p)\}}{\overline{F}(u_p)} \leq \varphi(1-\vep).
\end{align*}
Since $\vep > 0$ is arbitrary and $\varphi$ continuous at $1$ it is possible
to let $\vep \to 0$ and get the
upper bound 
\begin{align*}
   \limsup_{p \to 1}\frac{E_{\nu_{u_p}}w_{u_p}(X_1)^2
      I\{X_1 > F^\leftarrow(p)\}}{\overline{F}(u_p)} \leq \varphi(1).
\end{align*}

Then, by \eqref{eq:as.rnscale2}, and
\eqref{eq:as.rnscale3} it follows that
\begin{align*}
  \limsup_{p \to 1}\frac{\sigma^2_p}{F^\leftarrow(p)^2} &=
   \limsup_{p \to 1} \frac{E_{\nu_{u_p}}w_{u_p}(X_1)^2
    I\{X_1 > F^{\leftarrow}(p)\} - (1-p)^2}{F^\leftarrow(p)^2
    f(F^\leftarrow(p))^2} \\ 
  &=  \limsup_{p \to 1} \frac{E_{\nu_{u_p}}w_{u_p}(X_1)^2
    I\{X_1 > F^{\leftarrow}(p)\}}{\overline{F}(U^{\leftarrow}(p))^2}
  \frac{\overline{F}(U^{\leftarrow}(p))^2}{F^\leftarrow(p)^2
    f(F^\leftarrow(p))^2} \\
  & \quad - \frac{(1-p)^2}{F^\leftarrow(p)^2
    f(F^\leftarrow(p))^2} =  \frac{\varphi(1)-1}{\alpha^2}.
\end{align*}
\end{proof}

\subsection{Expected Shortfall}
Next we consider the properties, when $p$ is close to $1$, of the asymptotic
variance in the central limit theorem, Proposition \ref{prop:es}, for expected
shortfall.  
\begin{prop}\label{prop:es1}
  Let $\alpha > 2$.   Suppose \eqref{eq:as.regvar} and
  \eqref{eq:as.approx} hold. Suppose additionally that 
  there exist $c_0 < 1$ and a non-increasing function $\varphi$, regularly
  varying with index $-\alpha$, such that, for $c \geq c_0$, 
  \begin{align}\label{eq:es.rv}
    \limsup_{\lambda \to \infty}\frac{E_{\nu_\lambda}w_\lambda(X_1)^2 
      I\{X_1 > c\lambda\}}{\overline{F}(\lambda)^2} \leq \varphi(c). 
  \end{align}
  Then the sampling measures $\nu_{u_p}$ satisfy
  \begin{align*}
    \limsup_{p \to 1} \frac{\Var\Big(\frac{1}{1-p}\int_p^1
      \frac{Z(F^{\leftarrow}(u))}{f({F}^{\leftarrow}(u))}du\Big)}{\gamma_p(F^\leftarrow)^2}       
    < \infty.  
  \end{align*}
  Moreover, if \eqref{eq:es.rv} holds with $h(c) = K c^{-\alpha}$, for some
  constant $K \in (0,\infty)$, then
  \begin{align*}
    \limsup_{p \to 1} \frac{\Var\Big(\frac{1}{1-p}\int_p^1
      \frac{Z(F^{\leftarrow}(u))}{f({F}^{\leftarrow}(u))}du\Big)}{\gamma_p(F^\leftarrow)^2}
    \leq \frac{1}{\alpha^2}\Big(\frac{2K(\alpha-1)}{\alpha-2} - 1\Big).
  \end{align*}
\end{prop}

\begin{proof}
 By Lemma \ref{lem:covfcn},
\begin{align*}
  &\varrho(F^{\leftarrow}(u),F^{\leftarrow}(v))\\
  &\quad = E_{\nu_{u_p}}\Big( w_{u_p}(X_1)^2 I\{
  X_1\!>\!F^{\leftarrow}(u) \}I\{X_1\!>\!F^{\leftarrow}(v)\}\Big)\!-\!\overline{F}(F^{\leftarrow}(u))\overline{F}(F^{\leftarrow}(v))\\   
  &\quad = E_{\nu_{u_p}}\Big( w_{u_p}(X_1)^2 I\{
  X_1>F^{\leftarrow}(u) \}I\{X_1>F^{\leftarrow}(v) \}\Big)-(1-u)(1-v),
\end{align*}  
which implies, by Proposition \ref{prop:tedfconv}, that
\begin{align}
  & \frac{\Var\Big( \frac{1}{1-p} \int_p^{1}
    \frac{Z(F^{\leftarrow}(q))}{f(F^{\leftarrow}(q))} 
  dq\Big)}{\gamma_p(F^\leftarrow)^2} \nonumber \\ &=
\frac{1}{(1-p)^2\gamma_p(F^\leftarrow)^2}  \int_p^{1}\int_p^{1} 
  \frac{E(Z(F^{\leftarrow}(u))Z(F^{\leftarrow}(v)))}{f(F^{\leftarrow}(u))f(F^{\leftarrow}(v))} 
  dudv \nonumber \\ 
   &= \frac{1}{(1-p)^2\gamma_p(F^\leftarrow)^2}  \int_p^{1}\int_p^{1}
  \frac{\varrho(F^{\leftarrow}(u),F^{\leftarrow}(v))}{f(F^{\leftarrow}(u))f(F^{\leftarrow}(v))}
  dudv \nonumber \\
&= \frac{1}{(1-p)^2\gamma_p(F^\leftarrow)^2} \nonumber  \\ & \qquad \times \int_p^{1}\int_p^{1}
  \frac{E_{\nu_{u_p}}\Big( w_{u_p}(X_1)^2 I\{
  X_1>F^{\leftarrow}(u) \}I\{X_1>F^{\leftarrow}(v)
  \}\Big)}{f(F^{\leftarrow}(u))f(F^{\leftarrow}(v))} 
  dudv \label{eq:line1}\\
  & \quad - \frac{1}{(1-p)^2\gamma_p(F^\leftarrow)^2}  \int_p^{1}\int_p^{1} 
  \frac{(1-u)(1-v)}{f(F^{\leftarrow}(u))f(F^{\leftarrow}(v))} 
  dudv.\label{eq:line2}
\end{align}
Consider first \eqref{eq:line2}. By Karamata's theorem $f(F^\leftarrow(u)) \sim \frac{\alpha}{
F^\leftarrow(u)}\overline{F}(F^\leftarrow(u)) = \frac{\alpha(1-u)}{
F^\leftarrow(u)}$, as $u \to 1$, and therefore  
\begin{align*}
	\lim_{p\to 1}
        \frac{\frac{1}{(1-p)}\int_p^1\frac{(1-u)}{f(F^\leftarrow(u))}
          du}{\gamma_p(F^\leftarrow)} = 	\lim_{p\to 1}
        \frac{\int_p^1 \alpha^{-1}
F^\leftarrow(u)du}{\int_p^1 F^\leftarrow(u)du} = \alpha^{-1}. 
\end{align*}
This yields,
\begin{align*}
  \lim_{p\to 1}
  \frac{\frac{1}{(1-p)^2}\int_p^1\int_p^1
    \frac{(1-u)(1-v)}{f(F^\leftarrow(u))f(F^\leftarrow(v))} 
    dudv}{\gamma_p(F^\leftarrow)^2} = \alpha^{-2}.
\end{align*}
Next rewrite \eqref{eq:line1} for $v>u$ as 
\begin{align*}
  & \frac{2}{(1-p)^2\gamma_p(F^\leftarrow)^2}\\ & \quad \times \int_p^{1} 
  \frac{1}{f(F^{\leftarrow}(u))}\Big[\int_u^{1} 
  \frac{E_{\nu_{u_p}}\Big( w_{u_p}(X_1)^2I\{X_1>F^{\leftarrow}(v)
  \}\Big)}{f(F^{\leftarrow}(v))} 
  dv\Big]du 
\end{align*}
Then, the inner integrand can be written 
\begin{align}
  &\frac{E_{\nu_{u_p}}\Big( w_{u_p}(X_1)^2I\{X_1>F^{\leftarrow}(v)
  \}\Big)}{f(F^{\leftarrow}(v))} \nonumber \\
&\quad = \frac{E_{\nu_{u_p}}\Big( w_{u_p}(X_1)^2I\{X_1>\frac{F^{\leftarrow}(v)}{u_p}u_p\}\Big)}{\overline{F}(u_p)^2} 
\frac{\overline{F}(u_p)^2}{\overline{F}(F^\leftarrow(p))^2}
\frac{\overline{F}(F^\leftarrow(p))^2}{f(F^{\leftarrow}(v))} \nonumber
\\
& \quad \lesssim
\varphi\Big(\frac{F^\leftarrow(v)}{u_p}\Big)\frac{(1-p)^2}{f(F^\leftarrow(v))},  
\label{eq:es.inner}
\end{align}
where we have used \eqref{eq:es.rv}. By Potter's bound there
exists, for each 
$\vep > 0$, a constant $C_\vep$ such that
$\varphi(F^\leftarrow(v)/u_p) \leq C_\vep
(F^\leftarrow(v)/u_p)^{-\alpha+\vep}$. Take $ 0 < \vep <
2-\alpha$. The asymptotics of the integral in
\eqref{eq:line1} can now 
be determined as 
\begin{align*}
  &\int_p^{1} 
  \frac{1}{f(F^{\leftarrow}(u))}\Big[\int_u^{1} 
  \frac{E_{\nu_{u_p}}\Big( w_{u_p}(X_1)^2I\{X_1>F^{\leftarrow}(v)
  \}\Big)}{f(F^{\leftarrow}(v))} 
  dv\Big]du \\
  & \quad \lesssim  
  C_\vep u_p^{\alpha-\vep} (1-p)^2 \int_p^1
  \frac{1}{f(F^\leftarrow(u))} \int_u^1
  \frac{F^\leftarrow(v)^{-\alpha+\vep}}{f(F^\leftarrow(v))}dv \\
  & \quad = C_\vep u_p^{\alpha-\vep} (1-p)^2 \int_p^1
  \frac{1}{f(F^\leftarrow(u))} \int_{F^\leftarrow(u)}^\infty
  y^{-\alpha+\vep} dydu \\
  & \quad = C_\vep u_p^{\alpha-\vep} (1-p)^2 \int_{F^\leftarrow(p)}^\infty
  \frac{x^{1-\alpha+\vep}}{\alpha-\vep-1} dx \\
  & \quad = C_\vep u_p^{\alpha-\vep} (1-p)^2
  \frac{F^\leftarrow(p)^{2-\alpha+\vep}}{(\alpha-\vep-1)(\alpha-\vep-2)}\\
  & \quad \sim C_\vep (1-p)^2
  \frac{F^\leftarrow(p)^{2}}{(\alpha-\vep-1)(\alpha-\vep-2)}.
\end{align*}
By Karamata's theorem
\begin{align*}
\gamma_p(F^\leftarrow) & = \frac{1}{1-p}\int_p^1 F^\leftarrow(u) du =
\frac{1}{1-p}\int_{F^\leftarrow(p)}^\infty x f(x) dx \\ & \sim
\frac{1}{1-p}\frac{\alpha}{\alpha-1}F^\leftarrow(p)
  \overline{F}(F^\leftarrow(p)) \\
& = \frac{\alpha}{\alpha-1}F^\leftarrow(p).
\end{align*}
Putting everything together, the expression in \eqref{eq:line1} is
asymptotically bounded.

Moreover, if \eqref{eq:es.rv} holds with $\varphi(c) = K c^{-\alpha}$, for some
constant $K \in (0,\infty)$, then it is possible to take $C_\vep = K$
and $\vep = 0$. This results in
\begin{align*}
  \limsup_{p \to 1} \frac{\Var\Big(\frac{1}{1-p}\int_p^1
    \frac{Z(F^{\leftarrow}(u))}{f({F}^{\leftarrow}(u))}du\Big)}{\gamma_p(F^\leftarrow)^2}
  \leq \frac{1}{\alpha^2}\Big(\frac{2K(\alpha-1)}{\alpha-2} - 1\Big).
\end{align*}
\end{proof}

\section{Examples and numerical illustrations}
\label{sec:num}

In this section we use the methods presented in the previous
sections to design efficient algorithms for computing Value-at-Risk
and expected shortfall of a random variable $X$ which is the value
at time $n \geq 1$ of a heavy-tailed random walk. More precisely, 
\begin{align}\label{eq:x}
  X = \sum_{i=1}^n Z_i,
\end{align}
where $Z_i, \, i=1,\ldots , n$, are i.i.d.\ and regularly varying with
tail index $\alpha$. We will use $F_X$ and $F_Z$ to denote the d.f.\
of $X$ and $Z_1$, respectively. We write $f_X$ and $f_z$ for the
corresponding densities. First we need to establish that the
assumptions in the beginning of Section 4 
are satisfied. 

The subexponential property implies that the tail 
the random variable $X$ satisfies $\overline{F}_X(x) \sim n
\overline{F}_Z(x)$, as $x 
\to \infty$. Hence, $\overline{F}_X$ is regularly varying with index
$-\alpha$ and the function $U$ can be taken to be $n
\overline{F}_Z$. Finally, we need to consider importance sampling
algorithms with 
bounded relative error for computing rare event probabilities of the
form $P(X > \lambda)$. 

There exist several importance sampling
algorithms for efficient computation of rare event probabilities of
this form. Here we consider the dynamic mixture algorithms described
in \cite{HS09a} to generate $N$ independent samples of $X$,
denoted $X_1, \dots, X_N$. In particular we consider the conditional
mixture algorithm of \cite{DLW07} and the scaling mixture
algorithm of \cite{HS09a}.  Then, the tail e.d.f.\
$\overline{F}_{\nu_{u_p},N}$ is constructed from the sample. 
Value-at-Risk is computed as 
$(1-\overline{F}_{\nu_{u_p},N})^\leftarrow(p)$ and expected shortfall
as $\gamma_p((1-\overline{F}_{\nu_{u_p,N}})^\leftarrow)$. 

In the next subsection we verify the conditions of Proposition
\ref{prop:es}, Proposition 
\ref{prop2:rnscale}, and Proposition \ref{prop:es1} for these
algorithms. Then the algorithms are implemented and their numerical
performance is illustrated when $Z_1$ has a Pareto distribution. 

\subsection{Dynamic mixture algorithms}
The dynamic mixture algorithm is designed for generating samples of
$X$ in \eqref{eq:x} in order to efficiently compute rare event probabilities of
the form $P(X > \lambda)$. Here 
it is convenient to use the notation $S_i = Z_1 + \dots + Z_i$, $i
\geq 1$, $S_0 = 0$, and with this notation $X = S_n$ is the variable
of interest. Each sample of $S_n$ is generated sequentially by
sampling $Z_i$ from a mixture where the distribution of $Z_i$ may
depend on the current state, $S_{i-1}$. In the $i$th step,
$i=1,\dots,n-1$, where $S_{i-1} = s_{i-1}$,  
$Z_i$ is sampled as follows.
\begin{itemize}
\item If $s_{i-1} > \lambda$, $Z_i$ is sampled from the original
  density $f_Z$,
\item if $s_{i-1} \leq \lambda$, $Z_i$ is sampled from
  \begin{align*}
    &p_if_Z(\cdot) + q_i g_i( \cdot \mid s_{i-1}), \quad \text{ for }1
    \leq i \leq n-1,\\ 
    &g_n(\cdot \mid s_{n-1}), \quad \text{ for } i = n,
\end{align*}
where $g_i( \cdot \mid s_{i-1})$ is a state dependent density. Here
$p_i + q_i = 1$ and $p_i \in (0,1)$.  
\end{itemize}
The sampling measure distribution of $S_n$ obtained by the dynamic
mixture algorithm for computing $P(S_n > \lambda)$ is, throughout this
section, denoted $\nu_\lambda$. 

The following results provide sufficient conditions for the upper bound
$\varphi(c)$ that appears in Proposition \ref{prop2:rnscale}
and Proposition \ref{prop:es1}. 

\begin{lem}\label{lem:RNbddgen} 
  Consider the mixture algorithm above with $p_i > 0$ for $1 \leq i
  \leq n-1$. Suppose there exist $a \in (0,1)$ and $c>0$ such that
  \begin{align}
    \liminf_{\lambda \to \infty}\inf_{\tiny \begin{array}{l} s \leq
        c(1-(1-a)^{i})\\ y > a(c-s)\end{array}}\frac{g_i(\lambda y\mid
      \lambda s)}{f_Z(\lambda y)} 
    \overline{F}_Z(\lambda) &> 0, \quad  1 \leq i \leq n, \label{eq:main}\\
    \limsup_{\lambda \to \infty}\sup_{\tiny \begin{array}{l} s \leq c\\ y >
        c-s \end{array}}\frac{f_Z(\lambda y)}{g_n(\lambda y\mid \lambda
      s)} &< \infty.\label{eq:sec} 
  \end{align}
  Then the scaled 
  Radon-Nikodym derivative
  $\frac{1}{\overline{F}_Z(\lambda)}\frac{d\mu}{d\nu_\lambda}(\lambda y)$ is
  bounded on $\{y_1+\dots + y_n > c\}$.    
\end{lem}
The proof is essentially identical to the proof of Lemma 3.1 in
\cite{HS09a} and therefore omitted. 

\begin{thm}\label{thm:main}
  Suppose \eqref{eq:main} and \eqref{eq:sec} hold for $a \in
  (0,1)$. Suppose, in addition, that
  there exist continuous functions $h_i: \R^n \to [0,\infty)$ and a
  constant $c_0>0$ such that 
  \begin{align}\label{eq:thm}
    \frac{f_Z(\lambda y_i)}{g_i(\lambda y_i \mid \lambda s_{i-1})
      \overline{F}(\lambda)} \to h_i(y_i
    \mid s_{i-1}), 
  \end{align}
  uniformly on $\{y \in \R^n: s_{i-1} \leq c(1-(1-a)^{i-1}), y_i >
  a(c-s_{i-1})\}$ for any $c\geq c_0$.   
Then, for $c \geq c_0$, 
\begin{align}\label{eq:phi}
  \limsup_{\lambda \to\infty} \frac{E_{\nu_\lambda}w_\lambda(X)^2I\{X
    > c\lambda\}}{\overline{F}_X(\lambda)^2} \leq 
  \sum_{i=1}^n \prod_{j=1}^{i-1}\frac{1}{p_j} \frac{1}{q_i}
  \int_c^\infty  h_i(y_i \mid 0) \alpha y_i^{-\alpha-1}dy_i,
\end{align}
with $q_n = 1$. 
\end{thm}
The proof is essentially identical to the proof of Theorem 3.2 in
\cite{HS09a} and therefore omitted. 

By Theorem \ref{thm:main} we
see that the function $\varphi$ in Proposition \ref{prop2:rnscale} can
be taken as the right-hand side of \eqref{eq:phi}. Moreover, by
Karamata's theorem, it is regularly varying with index 
$-\alpha$ if $h_i(y_i\mid 0)$ is slowly varying. 

Finally, we establish that the conditions on the covariance function
$\varrho$ in Proposition \ref{prop:es} are satisfied.  
\begin{lem}
  Let $X$ have distribution $\mu$, d.f.\ $F_X$ and density $f_X$. 
  Suppose  $\overline{F}_X$ is regularly varying with index $-\alpha$, with
  $\alpha > 2$. Let $\nu$ denote any sampling distribution.
  If $d\mu/d\nu$ is bounded on $(a,\infty)$ then the covariance
  function $\varrho$ in \eqref{eq:rho} satisfies
  $\int_a^\infty \int_a^\infty \varrho(x,y)dxdy < \infty$ and
  $\varrho(x,x) = o([f_X(x)/\overline{F}_X(x)]^2)$. 
\end{lem}
\begin{proof}
  First note that if $d\mu/d\nu \leq K$ for some constant $K \in
  (0,\infty)$ then 
  \begin{align*}
    \varrho(x,y) \leq K \overline{F}_X(y)-
    \overline{F}_X(x)\overline{F}_X(y), \quad \text{ for } y \geq x.  
  \end{align*}
  Then
  \begin{align*}
    \int_a^\infty \int_a^\infty \varrho(x,y)  dxdy &= 2
    \int_a^\infty \int_x^\infty \varrho(x,y) dy dx \\
    & \leq 2 K \int_a^\infty \int_x^\infty \overline{F}_X(y) dy dx - 2
    \Big(\int_a^\infty\overline{F}_X(y) dy\Big)^2. 
  \end{align*}
  The first integral is finite, by Karamata's theorem, since $\alpha >
  2$ and the second 
  integral is finite for $\alpha > 1$ and then also for $\alpha > 2$. 
  For the second condition  
  \begin{align*}
    \frac{\varrho(x,x)}{f_X(x)^2/\overline{F}_X(x)^2} \leq K
    \frac{\overline{F}_X(x)^3}{f_X(x)^2}-
    \frac{\overline{F}_X(x)^4}{f_X(x)^2}. 
  \end{align*}
  By Karamata's theorem $\alpha \overline{F}_X(x) \sim x f_X(x)$ so
  the expression in the last display is asymptotically equivalent to
  \begin{align*}
    K x^{3}f_X(x) - x^4 f_X(x).
  \end{align*}
  This converges to $0$ as $x \to \infty$ when $\alpha > 2$ since $f_X$
  is regularly varying with index $-\alpha-1$. This
  completes the proof.
\end{proof}

\subsection{Conditional mixture algorithms}
The conditional mixture algorithm by \cite{DLW07} can be treated
with the above results. 


The conditional mixture algorithm  has, with $a \in (0,1)$,
\begin{align*}
  g_i(x \mid s) &= \frac{f_Z(x)I\{x > a(b-s)\}}{\overline{F}_Z(a(b-s))},
  \quad 1 \leq i \leq n-1,\\
  g_n(x \mid s) &= \frac{f_Z(x)I\{x > b-s\}}{\overline{F}_Z(b-s)}.
\end{align*}
Then the techniques for establishing the conditions of Lemma
\ref{lem:RNbddgen} and Theorem \ref{thm:main}, with $c_0 = 1$, are
completely similar to the ones in Section 4.1 in
\cite{HS09a}. The upper bound in Theorem \ref{thm:main} holds where
the functions $h_i$ are given by (see \cite{HS09a})
\begin{align*}
  h_i(y\mid s) &=
  \lim_{\lambda \to\infty} \frac{f_Z(by)}{f_Z(\lambda y)/\overline{F}_Z(a\lambda(1-s))
    \overline{F}_Z(\lambda)I\{y>a(1-s)\}}\\
  &=a^{-\alpha}(1-s)^{-\alpha}, \quad i=1,\ldots n-1
\end{align*}
and
\begin{align*}
  h_n(y\mid s)=\lim_{\lambda \to\infty} \frac{f_Z(\lambda
    y)}{f_Z(\lambda y)/\overline{F}_Z(\lambda(1-s))
    \overline{F}_Z(\lambda)I\{y>(1-s)\}}=(1-s)^{-\alpha}.
\end{align*}
The resulting upper bound $\varphi(c)$ in Propositions
\ref{prop2:rnscale} and \ref{prop:es1} is given by  
\begin{align*}
  \varphi(c)=c^{-\alpha}\bigg(a^{-\alpha}\sum_{i=1}^{n-1} \prod_{j=1}^{i-1}\frac{1}{p_j} 
  \frac{1}{q_i}+\prod_{j=1}^{n-1}\frac{1}{p_j}\bigg),
\end{align*}
where $a\in (0,1)$.
  
\subsection{Scaling mixture algorithms}
In the scaling mixture algorithm  the large variables are generated by
sampling from the original density and multiplying
with a large number. In the context of scaling mixtures we assume that
the orginal density $f_Z$ is strictly positive on $(0,\infty)$. 
We also assume that $f_Z(x) = x^{-\alpha-1}L(x)$ with $L$
slowly varying and $\inf_{x> x_0}L(x) =: L_* > 0$ for some $x_0>0$. 
The scaling mixture algorithm, with $\sigma > 0$, has
\begin{align*}
  g_i(x\mid s) &= (\sigma \lambda)^{-1}f_Z(x/\sigma \lambda)I\{x > 0\} + f_Z(x)I\{x
  \leq 0\},\quad i=1,\dots,n-1,\\
  g_n(x \mid s) &=  (\sigma \lambda)^{-1}f_Z(x/\sigma \lambda)I\{x > 0, s \leq
  \lambda-\lambda(1-a)^{n-1}\}\\  
  & \quad + f_Z(x)I\{x \leq 0 \text{ or } s > \lambda-\lambda(1-a)^{n-1}\}.
\end{align*}
To generate a sample $Z$ from $g_i$ proceed as follows. Generate a candidate
$Z'$ from $f_Z$. If $Z' \leq 0$ put $Z = Z'$ and if $Z' > 0$, put $Z =
\sigma \lambda Z'$. 

For the scaling mixture algorithm the conditions of Lemma
\ref{lem:RNbddgen} and Theorem \ref{thm:main} can be established with
$c_0 < 1$. The techniques for doing this are
completely similar to the ones in Section 4.3 in
\cite{HS09a}. The upper bound in Theorem \ref{thm:main} holds where
the functions $h_i$ are given by (see \cite{HS09a})
\begin{align*}
  h_i(y_i \mid s_{i-1}) = \alpha \lambda [y_i^{\alpha+1}f(y_i/\lambda)]^{-1}.
\end{align*}

  and the resulting upper bound $\varphi(c)$ in Propositions
  \ref{prop2:rnscale} and \ref{prop:es1}  is given by
\begin{align*}
  \varphi(c)=\sum_{i=1}^n \prod_{j=1}^{i-1}\frac{1}{p_j} \frac{1}{q_i}
  \int_c^\infty   \frac{\alpha }{\lambda^{\alpha}
    L(y_i/\lambda)}\alpha y_i^{-\alpha-1}dy_i. 
\end{align*}

\subsection{Numerical computation of Value-at-risk}
We now consider a sum of $n$ Pareto-distributed random variables, 
\begin{align*}
  S_{n}=Z_1+\ldots +Z_{n}.
\end{align*}
We will estimate quantiles of $S_n$ by using the importance sampling e.d.f.~given by 
the scaling mixture algorithm in \cite{HS09a} (SM) as well as 
the conditional mixture algorithm in \cite{DLW07} (DLW). The changes of measure are 
chosen by using the asymptotic approximation of the quantiles,
\begin{align*}
  \lambda_p^*=(n/(1-p))^{1/\alpha}-1.
\end{align*}
This approximation is based on the subexponential property, and since $P(S_n>x)>nP(X_1>x)$ 
for positive random variables, it is smaller than the true quantile.

For $p$ equal to $0.99, 0.999$ and $0.99999$, we use the DLW algorithm $10^2$ times
with $N=5\cdot 10^4$ samples to obtain a reference value which we refer to 
as the true value of the quantile. 

We compare the performance of the quantile estimates based on $N=10^4$ 
samples. The estimation is repeated $100$ times and the mean and standard 
deviation of the estimates are reported.

We also include the results from standard Monte Carlo for comparison.


\begin{table}[ht]
\caption{Simulations of $\lambda_p$ such that $P(S_{n}>\lambda_p)=1-p$, 
where $S_n=\sum_{i=1}^n Z_i$ and $P(Z_1>x)=(1+x)^{-2}$.
The number of samples used for each estimate was $N=10^4$ and the estimation was repeated 100 times.}\label{table1}
 
\begin{tabular}{|c|c|c|c|c|c|c|c|}
\hline
    $n$ &  $1-p$ 	&   True 	& Approx.	& SM   		&  DLW		& MC 		&\\
\hline
    10	&  1e-2 	&   40.141 	&  30.623	&   41.007	& 40.166 	& 40.038 	& Avg.~est.\\
	&   		& 		&		&   (0.246)	& (0.459) 	& (1.780)	&(Std.~dev.)\\
\hline 
  	&  1e-3		&   108.49 	& 99.000	&  109.33  	& 108.29	&84.821		&\\
	&  		&    	 	&		&  (0.847) 	& (1.081)	&(47.23)	&\\
\hline
   	&  1e-5 	&   1007.4 	& 999.00	&  1003.1 	& 1007.5	& 609.42 	&\\
	&   		&  	 	&		&  (18.5)  	& (1.51)	& (1594)	&\\
\hline
  30 	&  1e-2		&  84.622	& 53.772 	&  85.841 	& 84.681	& 84.362	&\\
	&   		&   		&		&  (0.3950)  	& (1.237)	&  (2.739)	&\\
\hline
	&  1e-3		&  202.41	& 172.21	&   203.56  	& 202.29	&  171.16	&\\
	&  		&   		&		&   (1.530)	& (2.400)	&  (71.26) 	&\\
\hline
	&  1e-5		&  1759.5	&  1731.1	&  1753.7   	& 1759.0 	& 114.23	&\\
	&  		&   	 	&		&  (41.12)	& (1.487) 	& (443.5)	&\\
\hline

\end{tabular} 
\end{table}

\begin{table}[ht]
\caption{Simulations of $\lambda_p$ such that $P(S_{n}>\lambda_p)=1-p$, 
where $S_n=\sum_{i=1}^n Z_i$ and $P(Z_1>x)=(1+x)^{-3}$.
The number of samples used for each estimate was $N=10^4$ and the estimation was repeated 100 times.}\label{table2}
 
\begin{tabular}{|c|c|c|c|c|c|c|c|}
\hline
$n$ 	&  $1-p$ 	&   True 	& Approx.	& SM	   	&  DLW		& MC 		&\\
\hline
    10	&  1e-2 	&    14.190 	& 9.0000	&  14.853	& 14.195 	& 14.182 	& Avg.~est.\\
	&   		& 		&		&  (0.090)	& (0.154) 	& (0.305)	&(Std.~dev.)\\
\hline 
  	&  1e-3		&    25.656	&  20.544	&  26.125  	&  25.588	&24.965		&\\
	&  		&    	 	&		& (0.171) 	& (0.412)	&(2.212)	&\\
\hline
   	&  1e-5 	&   103.42	&  99.000	& 104.23	& 103.40	& 5.283 	&\\
	&   		&  	 	&		& (0.799)  	& (0.553)	& (16.03)	&\\
\hline
  30 	&  1e-2		&   29.951	&  13.422 	& 31.054 	& 29.943	& 29.949	&\\
	&   		&   		&		& (0.287)  	& (0.519)	&(0.500)	&\\
\hline
	&  1e-3		&   46.072	& 30.072	&  46.725  	& 46.277	&  44.608	&\\
	&  		&   		&		&   (0.286)	& (1.041)	&  (2.688) 	&\\
\hline
	&  1e-5		&   157.65	&  143.22	&  158.46  	& 157.62 	& 13.847	&\\
	&  		&   	 	&		&   (1.080)	& (0.273) 	& (28.53)	&\\
\hline

\end{tabular} 
\end{table}

\subsection{Numerical computation of expected shortfall}
Using the setting from the previous section, we also calculate the
expected shortfall for the case of a random walk with
Pareto-distributed increments. We first consider the case where $\alpha=2$, although 
it does not satisfy the conditions of Proposition \ref{prop:es1}.

\begin{table}[ht]
\caption{Simulations of $E(S_n|S_n>\lambda_p$), where $P(S_{n}>\lambda_p)=1-p$, 
$S_n=\sum_{i=1}^n Z_i$ and $P(Z_1>x)=(1+x)^{-2}$. 
The number of samples used for each estimate was $N=10^4$ and the estimation was repeated 100 times.}\label{table3}
 
\begin{tabular}{|c|c|c|c|c|c|c|}
\hline
$n$ 	&  $1-p$ 	&    True value	& SM		& DLW 		&  MC   	&\\
\hline
      10&  1e-2 	&   71.795	&  73.065	& 71.831	&   72.252	& Avg.~est.\\
	&   		& 		&  (1.06)	& (1.22)	&   (8.75)	&(Std.~dev.)\\
	&		&		&  [0.845]	& [0.815]	&  [0.702]	& [Avg.~time (s)]	\\
\hline 
	  &  1e-3	&   208.84	&  209.37	& 209.30	&   213.42	&\\
	&  		&   		& (3.60) 	& (4.99)	&  (65.8)	&\\
	&		&		& [0.734]	& [0.724]	&  [0.572]	&\\
\hline
   	&  1e-5		&    2008.4	&   2009.8  	&   2009.3	& 4787.8	&\\
	&	   	&  		&  (37.1)	&   (30.9)	&  (23168) 	&\\
	&		&		&  [0.866]	&   [0.822]	&  [0.693]	&\\
\hline
  30 	&  1e-2		&   139.22	&   140.55  	&  139.14	&  140.76	&\\
	&	   	&  		&  (2.22)	&   (3.09)	&  (17.34) 	&\\
	&		&		&  [1.189]	&   [1.077]	&  [0.903]	&\\
\hline
	&  1e-3		&   376.29	&   375.76	&   378.24	&  391.06    	&\\
	&  		&  		&   (5.00)	&   (11.49)	& (96.36)	&\\
	&		&		&  [1.033]	&   [0.936]	& [0.757]	&\\
\hline
	&  1e-5		&   3494.4 	&  3500.0	&  3496.9	&  745.70 	&\\
	&  		&   	 	&  (65.2)	&   (59.8)	&  (3671) 	&\\
	&		&		&  [1.301]  	&  [1.224]	&  [0.991]	&\\

\hline

\end{tabular} 
\end{table}

\begin{table}[ht]
\caption{Simulations of $E(S_n|S_n>\lambda_p$), where $P(S_{n}>\lambda_p)=1-p$, 
$S_n=\sum_{i=1}^n Z_i$ and $P(Z_1>x)=(1+x)^{-3}$. 
The number of samples used for each estimate was $N=10^4$ and the estimation was repeated 100 times.}\label{table4}
 
\begin{tabular}{|c|c|c|c|c|c|c|}
\hline
$n$ 	&  $1-p$ 	&   True value	& SM		& DLW 		&  MC   	&\\
\hline
      10&  1e-2 	&   19.260	&  20.044	&  19.257	&    19.495	& Avg.~est.\\
	&   		& 		&  (0.167)	&  (0.395)	&   (0.905)	& (Std.~dev.)\\
	&		&		&  [0.702]	&  [0.727]	&   [0.605]	& [Avg.~time (s)]\\
\hline 
	  &  1e-3	&   36.658	&  36.911	&  36.463	&   41.032	&\\
	&  		&   		&  (0.327) 	&  (0.776)	&  (5.53)	&\\
	&		&		&  [0.7079]	&  [0.731]	&  [0.613]	&\\
\hline
  	&  1e-5 	&   154.74  	&  154.39	&  153.83	&   132.74 	&\\
	&   		& 	 	&  (1.326)	&  (2.705)	&   (491.7)  	&\\
	&		&		&  [0.708]	&  [0.733]	&   [0.607]	&\\
\hline
  30 	&  1e-2		&   37.277	&  38.603 	&  37.200	&   37.744	&\\
	&	   	&  		&  (0.902)	&  (1.169)	&   (1.581) 	&\\
	&		&		&  [0.885]	&  [0.923]	&   [0.707]	&\\
\hline
	&  1e-3		&   62.090	&   62.013	&  62.066	&  69.369    	&\\
	&  		&  		&  (0.416)	&  (1.814)	&  (7.973)	&\\
	&		&		&  [0.912]	&  [0.939]	&  [0.712]	&\\
\hline
	&  1e-5		&   232.01 	&  230.27	&  230.00 	&  225.14	&\\
	&  		& 	 	&  (1.92)     	&  (1.47)	&  (932) 	&\\
	&		&		&  [0.911]	&  [0.935]	&   [0.703]	&\\
\hline

\end{tabular} 
\end{table}

\end{document}